\documentclass[12pt]{article}
\usepackage{amssymb}
\usepackage{amsmath}

\newcommand\F{{\cal F}}

\newcommand\A{{\cal A}}

\newtheorem{thm}{Theorem}

\newtheorem{prob}{Open problem}
\newtheorem{obs}{Observation}
\newtheorem{cons}{Construction}

\title{Around the Complete Intersection Theorem}
\author{{\bf Gyula O.H. Katona}\thanks{The work of the author was
supported by the Hungarian National Foundation for Scientific
Research, grant number NK104183. }\\
MTA R\'enyi Institute\\
Budapest Pf 127, 1364 Hungary \\ ohkatona@renyi.hu }

\begin{document}
\date{}

\maketitle

\section{Introduction}

The underlying set will be $[n]=\{ 1,2,\ldots ,n\}$. The family of
all $k$-element subsets of $[n]$ is denoted by ${[n]\choose k}$. Its
subfamilies are called {\it uniform}. A family $\F$ of some subsets
of $[n]$ is called {\it intersecting} if $F\cap G\not= \emptyset$
holds for every pair $F,G \in \F$. It is easy to determine the
largest (non-uniform) intersecting family in $2^{[n]}$ since at most
one of the complementing pairs can be taken.

\begin{obs} {\rm(Erd\H os, Ko, Rado \cite{EKR})} If $\F \subset {2^{[n]}}$ is intersecting
then $|\F|\leq  2^{n-1}=2^n/2.$

\end{obs}

The following trivial construction shows that the bound is sharp.

\begin{cons} Take all subsets of $[n]$ containing the element 1.
\end{cons}

However there are many other construction giving equality in Observation 1. The following one
will be interesting for our further investigations.

\begin{cons} If $n$ is odd take all sets of size at least ${n+1\over 2}$.
If $n$ is even then choose
all the sets of size at least ${n\over 2}+1$ and the sets of size
${n\over 2}$  not containing the element $n$.
\end{cons}

The analogous problem when the intersecting subsets have size
exactly $k (k\leq {n\over 2})$, that is the case of uniform
families, is not so trivial.

\begin{thm}{\rm(Erd\H os, Ko, Rado \cite{EKR})} If $\F \subset {[n]\choose k}$ is intersecting
where $k\leq {n\over 2}$ then $$|\F|\leq  {n-1\choose k-1}.$$
\end{thm}

For a shorter proof see \cite{K2}. In this case there is only one extremal construction,
mimicking Construction 1.

\begin{cons} Take all subsets of $[n]$ having size $k$ and containing the element 1.
\end{cons}

Already Erd\H os, Ko and Rado, in their seminal paper considered a more general problem.
A family $\F \subset 2^{[n]}$ is $t${\it -intersecting} if
$|F\cap G|\geq t$ holds for every pair $F,G \in \F$. They posed a conjecture for the
maximal size of non-uniform
a $t$-intersecting family. This conjecture was justified in the following theorem.

\begin{thm}{\rm(Katona \cite{K1})} If $\F \subset 2^{[n]}$ is
$t$-intersecting  then
$$|\F|\leq \begin{cases} \sum_{i={n+t\over 2}}^n{n\choose i} &\textrm{ if } n+t
\textrm{ is even }\\ \sum_{i={n+t+1\over 2}}^n{n\choose
i}+{n-1\choose {n+t-1\over 2}}  &\textrm{ if } n+t \textrm{ is odd
}.\end{cases}$$
\end{thm}

Here the generalization of Construction 1 gives only $2^{n-t}$, less than the upper
bound in Theorem 2 (if $t>1$).
In order to obtain a sharp construction we have to mimic Construction 2.

\begin{cons} If $n+t$ is even, take all sets of size at least ${n+t\over 2}$.
If $n+t$ is odd then choose
all the sets of size at least ${n+t+1\over 2}$ and the sets of size ${n+t-1\over 2}$
not containing the element $n$.
\end{cons}

\section{The Complete Intersection Theorem}

The problem of $t$-intersecting families for the uniform case proved
to be much more difficult than for the non-uniform case. Erd\H os,
Ko and Rado were able to settle the problem when $n$ is large with
respect to $k$. The dependence of the threshold on $t$ is not
interesting here since $1\leq t< k$ can be supposed.

\begin{thm} {\rm  \cite{EKR}} If $\F \subset {[n]\choose k}$ is
$t$-intersecting and $n>n(k)$ then
$$|\F|\leq  {n-t\choose k-t}.\eqno(1)$$
\end{thm}

They also gave an example for small $n$ when (1) does not hold.
Let $n$ be divisible by 4, $k={n\over 2}$ and  $t=2$. The family
$$\F=\left\{ F:\ |F|={n\over 2}, \left| F\cap \left[ n\over 2\right] \right|
\geq \left[ n\over 4\right] +1 \right\}$$
is 2-intersecting, since any two members meet in at leats two
elements in $\left[ n\over 2\right]$. On the other hand the size of
this family is more than ${n-2\choose n/2-1}$, if $n>4$. (See e.g.
$n=8$.) They believed that this construction was optimal.

The next step towards the better understanding the situation was when Frankl \cite{F1}
and Wilson \cite{W} determined the exact value of the threshold $n(k)$ in Theorem 3.

Let us now consider the following generalization of the construction
above.
\begin{cons} Choose a non-negative integer parameter $i$ and
define the family
$$\A (n,k,t,i)=\{ A:\ |A|=k, | A\cap [t+2i]| \geq t+i\}. \eqno(2)$$
It is easy to see that $\A$ is $t$-intersecting for each $i$.
\end{cons}
Introduce the following notation:
$$\max_{0\leq i}|\A (n,k,t,i)|=\textrm{AK}(n,k,t).$$
This is the size of the best of the constructions (2). Frankl
\cite{F1} conjectured that this is construction gives the largest
$k$-uniform $t$-intersecting family. Frankl and F\"uredi \cite{FF}
proved the construction for a very large class of parameters but the
full conjecture remained open until 1996 when it became a theorem.

\begin{thm} {\rm ( {\bf The Complete Intersection Theorem, Ahlswede and Khachatrian}
\cite{AK})} Let $\F \subset {[n]\choose k}$ be a $t$-intersecting family. Then
$$|\F|\leq   \textrm{AK}(n,k,t)$$  holds.
\end{thm}

This theorem was a very important step in the progress of the Extremal Set Theory.
Its proof was a far-reaching generalization of the transformation method introduced
in \cite{EKR}. The author of the present paper must confess that he had mixed feeling
when he learned about the result. On the one hand he was happy that a new important
result of the theory came into life. On the other hand, however, he was a little
disappointed because he had the plan to solve the conjecture later when he had time
to devote all
his energies to the solution.

Of course Theorem 4 has many consequences. We will exhibit only one new result of us,
in Section 4,
where this theorem is used and plays a role even in the formulation of the statement.

\section{An Open Problem}

Even the best theorems do not stop the progress in science. In contrary, they raise
new questions. Let us show one.

If $\F \subset 2^{[n]}$ is a family of subsets, let $p_i(\F)$ denote the number
of $i$-element members of $\F$, that is, $p_i(\F)=\left| \F \cap {[n]\choose i} \right|$.
Then the vector $p(\F)=(p_0(\F), p_1(\F), \ldots ,p_n(\F) )$ is called the
{\it profile vector} of $\F$.

Take all profile vectors of $t$-intersecting families. They will
form a set of points with integer coordinates in the
$n+1$-dimensional Euclidean space. The vertices of the convex hull
of this set of points are called the {\it extreme points} of the
class of $t$-intersecting families. If some sets are deleted from a
$t$-intersecting family then the remaining family will also be
$t$-intersecting. Hence if $(p_0, p_1, \ldots ,p_n )$ is the profile
vector of a $t$-intersecting family and $q_i\leq p_i$ holds then
$(q_0, q_1, \ldots , q_n )$ is also a profile vector of a
$t$-intersecting family. An extreme point $(p_0, p_1, \ldots ,p_n )$
is called {\it essential} if  there is no other essential point
$(r_0, r_1, \ldots , r_n )$ satisfying $p_i\leq r_i$ for all $i$.
Let $E_n(t)$ denote the set of essential extreme points of the set
of profile vectors of all $t$-intersecting families.

It is easy to see that if $\alpha_j\geq 0$ are fixed constants then
$$\max_{\F \textrm{ is } t \textrm {-intersecting }} \sum_{i=0}^n\alpha_jp_j(\F)\eqno(3)$$
is attained for at least one essential extreme point. Therefore if
we want to determine the maximum in (3) it is sufficient to
calculate the linear combination of each of the vectors in $E_n(t)$
with the given $\alpha_j$s and find the largest one among these
values. Observe that if the coefficients are all zero except for a
fixed $k$ for which $\alpha_k=1$ then (3) gives the size of the
largest $k$-uniform $t$-intersecting family. On the other hand  if
the coefficients $\alpha_i=1 (0\leq i\leq n)$ are taken then (3)
gives the total number of sets in the family.

The essential extreme points were determined for the case $t=1$ in \cite{EFK} (Theorem 6).
(For an easier treatment of the theory see the paper of Gerbner \cite{G}.)
We have no place to give the full form
of the statement of this theorem. But it is easy to check that if $k\leq {n\over 2}$
then the largest $k$th coordinate in the essential extreme points is ${n-1\choose k-1}$
giving Theorem 1.  On the other hand, calculating the sums of the coordinates of the
essential extreme points we obtain the formula in Observation 1.

\begin{prob} Determine the essential extreme points of the $t$-intersecting families ($t>1$).
\end{prob}
Of course we know some of the extreme points. The one that maximizes
the linear combination $\sum_{i=0}^np_i(\F)$. It is determined by
Construction 4 for Theorem 2. This point is the ``farthest" one from
the origin. The difficulty lies in the determination of the extreme
points near the axes. Yet, the extreme points along the axes are
given by Theorem 4 and Construction 5.

\section{Union-intersecting families}

The following problem was asked by J\'anos K\"orner.

Let $\F\subset 2^{[n]}$ and suppose that if $F_1,F_2,G_1,G_2\in \F, F_1\not= F_2, G_1\not= G_2$
holds then
$$(F_1\cup F_2)\cap (G_1\cup G_2)\not=\emptyset.$$
What is the maximum size of such a family?

He conjectured that the following construction gives the largest one.

\begin{cons} If $n$ is odd then take all sets of size at least ${n-1\over 2}$.
If $n$ is even then choose all the sets of size at least ${n\over 2}$ and the
sets of size ${n\over 2}-1$  containing the element $1$.
\end{cons}

We solved the problem in a more general setting.
A family $\F\subset 2^{[n]}$ is called a {\it union-}$t${\it -intersecting}
if
$$|(F_1\cup F_2)\cap (G_1\cup G_2)|\geq t$$
holds for any four members such that $F_1\not= F_2, G_1\not= G_2$.

\begin{thm}{\rm (Katona-D.T. Nagy \cite{KN} )}
If $\F\subset 2^{[n]}$ is a union-$t$-intersecting family then
$$|\F|\leq \begin{cases} \sum_{i={n+t\over 2}-1}^n{n\choose i} &\textrm{ if } n+t
\textrm{ is even }\\
\sum_{i={n+t-1\over 2}}^n{n\choose i}+\textrm{AK}\left( n,{n+t-3\over 2},t\right)
&\textrm{ if } n+t \textrm{ is odd }.\end{cases}$$
\end{thm}

The following construction shows that the estimate is sharp.
\begin{cons} If $n+t$ is even, take all the sets with size at least
${n+t\over 2}-1$. Otherwise choose all the sets of size at least
${n+t-1\over 2}$ and the  sets of size ${n+t-3\over 2}$ following
Construction 5 where $k={n+t-1\over 2}$ and $i$ chosen to maximize
(3).
\end{cons}

Since the result contains the AK-function, it is obvious that Theorem 4 must be
used in the proof of this theorem.

As before, the uniform case is more difficult. Yet, we will treat it in an even more
general form. A family $\F\subset 2^{[n]}$ is called a $(u,v)$-{\it union-intersecting}
if for different members $F_1, \ldots ,F_u, G_1,\ldots ,G_v$ the following holds:
$$\left( \cup_{i=1}^uF_i\right) \cap \left( \cup_{j=1}^vG_j\right)\not= \emptyset .$$

\begin{thm} {\rm (Katona-D.T. Nagy \cite{KN})}
Let $1\leq u\leq v$ and suppose that the family $\F\subset {[n]\choose k}$ is a
$(u,v)$-union--intersecting family then
$$|\F|\leq {n-1\choose k-1}+u-1$$
holds if $n>n(k,v)$.
\end{thm}
The following construction shows that the estimate is sharp.
\begin{cons} Take all $k$-element subsets containing the element 1, and choose
$u-1$ distinct sets non containing 1.
\end{cons}
The theorem does not give a solution for small values.
\begin{prob} Is there an Ahlswede-Khachatrian type theorem here, too?
\end{prob}

\end{document}